\documentclass[11pt]{amsart}
\usepackage[T1]{fontenc}
\usepackage{lmodern}
\usepackage{amsmath,amssymb,amsthm,mathtools}
\usepackage[margin=1in]{geometry}
\usepackage{microtype}
\usepackage{needspace}
\usepackage[colorlinks=true,linkcolor=blue,citecolor=blue,urlcolor=blue]{hyperref}

\newtheorem{theorem}{Theorem}[section]
\newtheorem{proposition}[theorem]{Proposition}
\newtheorem{lemma}[theorem]{Lemma}
\newtheorem{corollary}[theorem]{Corollary}
\theoremstyle{definition}

\newcommand{\C}{\mathbb C}
\newcommand{\R}{\mathbb R}
\newcommand{\Q}{\mathbb Q}
\newcommand{\Z}{\mathbb Z}
\newcommand{\PP}{\mathbb P}
\newcommand{\G}{\boldsymbol\Gamma}
\newcommand{\E}{\mathcal E}
\DeclareMathOperator{\ord}{ord}
\DeclareMathOperator{\Res}{Res}
\allowdisplaybreaks[2]

\title[Generic finiteness of planar vortex relative equilibria]
{Generic finiteness of rotating relative equilibria in the planar $N$-vortex problem}
\author{Jinxin Xue}
\date{November 12, 2026}
\address{Rm 305, Wenbei Building, New Cornerstone Sciences Laboratory, Department of Mathematics Sciences, Tsinghua University, 100084}
\email{jxue@tsinghua.edu.cn}

\begin{document}
\begin{abstract}
For every $N\ge2$, we prove that Zariski-generic nonzero real circulations,
with any prescribed signs, admit at most $N^{2N}$ labeled, collision-free
rotating relative equilibria modulo translations, rotations in $SO(2)$,
and positive scaling. The excluded parameters lie in an explicit finite
union of quadratic hypersurfaces. 
The proof combines a logarithmic-residue obstruction to positive-dimensional
complex solution sets with B\'ezout degree bounds.
\end{abstract}
\maketitle

\section{Introduction}
We consider the following planar $N$-vortex problem. 
Fix $N\geq2$ and $\G=(\Gamma_1,\ldots,\Gamma_N)\in(\R^*)^N$.
A labeled, collision-free configuration $q=(q_1,\ldots,q_N)$ in $\R^2$
is a rotating relative equilibrium if there are $c\in\R^2$ and
$\lambda\in\R^*$ such that
\begin{equation}\label{eq:physical}
 \lambda(q_i-c)=F_i(q),\qquad
 F_i(q):=\sum_{j\ne i}\Gamma_j
       \frac{q_i-q_j}{|q_i-q_j|^2},\qquad 1\leq i\leq N.
\end{equation}
The sign convention is fixed by this equation. We identify configurations
under translations, rotations in $SO(2)$, and positive dilations. 

For $k=(k_{ij})_{1\le i<j\le N}\in\Z^{\binom N2}$, write
$\|k\|_\infty=\max_{i<j}|k_{ij}|$. Put
\begin{equation}\label{eq:exceptional}
\begin{split}
 D_N&=N^{2N},\\
 Q_k(\G)&=\sum_{i<j}k_{ij}\Gamma_i\Gamma_j,\\
 \E_N&=\bigcup_{\substack{k\in\Z^{\binom N2}\\0<\|k\|_\infty\leq2D_N}}
       \{\G\in\C^N:Q_k(\G)=0\}.
\end{split}
\end{equation}

\begin{theorem}\label{thm:signed}
The set $\E_N$ is a proper algebraic subset of $\C^N$, defined over $\Q$.
For every
\[
 \G\in(\R^*)^N\setminus\E_N,
\]
the set of labeled, collision-free rotating relative equilibria satisfying
\eqref{eq:physical}, modulo translations, rotations and positive scaling,
is finite and has at most $D_N$ elements. The conclusion also applies when
$\sum_i\Gamma_i=0$.

For each prescribed sign vector $\varepsilon\in\{1,-1\}^N$, the allowed
parameters form a nonempty Zariski-open intersection with
\[
 \mathcal O_\varepsilon
   =\{\G\in\R^N:\varepsilon_i\Gamma_i>0\text{ for every }i\}.
\]
This intersection is open, dense, and of full Lebesgue measure in
$\mathcal O_\varepsilon$.
\end{theorem}

\begin{corollary}[Zero total circulation]\label{cor:zero-total}
Let
\[
 H=\left\{\G\in\C^N:\sum_i\Gamma_i=0\right\}.
\]
Then $H\cap\E_N$ is a proper algebraic subset of $H$. Consequently, the
same bound $D_N$ holds for relatively Zariski-generic nonzero real
circulations on $H$. The good set is relatively open, dense, and of full
$(N-1)$-dimensional Lebesgue measure in every nonempty intersection
$H\cap\mathcal O_\varepsilon$.
\end{corollary}

The generic assumption is necessary. An example with a continuum of relative equilibria can be found in \cite{YuZhu2025}. 
\subsection{Literature review}
The finiteness problem for planar Newtonian central configurations is a longstanding
question in celestial mechanics and appears on Smale's list of problems for
the twenty-first century \cite[Problem~6]{Smale1998}.
A related finiteness problem arises for planar point vortices, whose
interaction Hamiltonian is logarithmic \cite[Introduction]{HamptonMoeckel2009}.
The logarithmic interaction comes
from the fundamental solution of the planar Laplacian. Logarithmic interactions also occur in
the plasma analogy for the Laughlin wave function in the fractional quantum
Hall effect \cite[Section 3.1.2]{Tong}.

Palmore announced almost-everywhere finiteness for positive circulations
\cite[Corollary~5.2]{Palmore1982}; his note does not supply a proof of the
generic nondegeneracy theorem used to obtain it. No part of the proof below
depends on that assertion. O'Neil's general-$N$ work treats static
equilibria, translations, and collinear relative equilibria
\cite{ONeil1987}; the scope is summarized in
\cite[Introduction]{HamptonMoeckel2009}. It should not be read as a theorem
about all planar rotating configurations with arbitrary signs.

For four vortices, Xiang Yu proves finiteness of rotating relative equilibria for
arbitrary nonzero real circulations \cite{Yu2023}. For five vortices,
Xiang Yu--Shuqiang Zhu obtain finiteness outside a codimension-two algebraic exceptional
set, and finiteness for every same-sign circulation vector
\cite[Theorems~1.2 and~1.4]{YuZhu2025}. O'Neil's more recent generic theorem
applies for arbitrary $N$ under a mirror-symmetry restriction, with at most
six vortices off the symmetry line, together with a special eight-vortex
case \cite{ONeil2025}. The point of the present
argument is the unrestricted number of vortices and the direct signed
extension of the full complex-fiber finiteness statement. 

\section*{Acknowledgment}
The author is supported by the National Natural Science Foundation of China (Grant No. 12271285, 11790273), and by the New Cornerstone Investigator program and Xiaomi Foundation. Moreover, the author is a member of IAS supported by the Shiing-Shen Chern Member Fund this semester. 

{\bf AI declaration: } The work is done with the help of ChatGPT 6.1. The author was inspired by \cite[Section 3.1.2]{Tong} and  suggested GPT the method of complexification which is specially adapted to the current problem and utilize the fact that the log potential is the plane Green function. GPT found the key argument, in particular the  residue argument, and also assisted in writing. 
\section{A normalization that retains both rotation signs}

Put $T=\sum_i\Gamma_i$, $L=\sum_{i<j}\Gamma_i\Gamma_j$ and
$r_{ij}=q_i-q_j$. Pairing the terms indexed by $(i,j)$ and $(j,i)$ gives
\begin{align}
 \sum_i\Gamma_iF_i
 &=\sum_{i<j}\Gamma_i\Gamma_j
       \frac{r_{ij}-r_{ij}}{|r_{ij}|^2}=0,\notag\\
 \sum_i\Gamma_iq_i\cdot F_i
 &=\sum_{i<j}\Gamma_i\Gamma_j
       \frac{(q_i-q_j)\cdot r_{ij}}{|r_{ij}|^2}=L.
 \label{eq:pair-identities}
\end{align}
For \eqref{eq:physical}, translating the rotation center to $0$ therefore gives
\begin{equation}\label{eq:center-virial}
 \sum_i\Gamma_i(q_i-c)=0,\qquad
 \lambda I_{\R}=L,\qquad
 I_{\R}=\sum_i\Gamma_i|q_i-c|^2.
\end{equation}
If $T\ne0$, then $c=T^{-1}\sum_i\Gamma_iq_i$; if $T=0$, the first
identity imposes $\sum_i\Gamma_iq_i=0$. The actual rotation data are unique:
\[
 \lambda(q_i-c)=\lambda'(q_i-c')\ (\forall i)
 \ \Longrightarrow\
 (\lambda-\lambda')(q_2-q_1)=0
 \ \Longrightarrow\ \lambda=\lambda',\quad c=c'.
\]

For fixed $\G\in(\C^*)^N$, define the complex fiber $X_{\G}$ by
\begin{equation}\label{eq:complex}
 z_i=\sum_{j\ne i}\frac{\Gamma_j}{w_i-w_j},\qquad
 w_i=\sum_{j\ne i}\frac{\Gamma_j}{z_i-z_j},\qquad 1\le i\le N,
\end{equation}
with
\begin{equation}\label{eq:gauge}
 z_2-z_1=1,\qquad
 s_{ij}:=(z_i-z_j)(w_i-w_j)\ne0
 \quad(1\le i<j\le N).
\end{equation}
Here $z_i,w_i$ are independent complex variables. The unsliced equations
are invariant under $(z,w)\mapsto(tz,t^{-1}w)$, $t\in\C^*$, and
$z_2-z_1=1$ meets each orbit exactly once.

\Needspace{7\baselineskip}
\begin{lemma}\label{lem:injection}
For $\G\in(\R^*)^N$, physical similarity classes satisfying
\eqref{eq:physical}, with either sign of $\lambda$, are in bijection with
\[
 \{(z,w)\in X_{\G}:w_i=b\overline{z_i}\ (\forall i),\quad
 b=w_2-w_1\in\R^*\}.
\]
\end{lemma}

\begin{proof}
Represent $q_i-c$ by $\zeta_i\in\C$. The force equations become
\[
 \lambda\zeta_i
 =\sum_{j\ne i}\frac{\Gamma_j}{\overline{\zeta_i}-\overline{\zeta_j}},
 \qquad
 \lambda\overline{\zeta_i}
 =\sum_{j\ne i}\frac{\Gamma_j}{\zeta_i-\zeta_j}.
\]
For $a=\zeta_2-\zeta_1\ne0$, set
\begin{equation}\label{eq:injection}
 z_i=\frac{\zeta_i}{a},\qquad
 w_i=\lambda a\overline{\zeta_i},\qquad
 b=\lambda|a|^2.
\end{equation}
Then
\begin{align*}
 \sum_{j\ne i}\frac{\Gamma_j}{w_i-w_j}
 &=\frac1{\lambda a}\sum_{j\ne i}
      \frac{\Gamma_j}{\overline{\zeta_i}-\overline{\zeta_j}}
   =\frac{\zeta_i}{a}=z_i,\\
 \sum_{j\ne i}\frac{\Gamma_j}{z_i-z_j}
 &=a\sum_{j\ne i}\frac{\Gamma_j}{\zeta_i-\zeta_j}
   =\lambda a\overline{\zeta_i}=w_i,
\end{align*}
and
\begin{equation}\label{eq:signed-distances}
 z_2-z_1=1,\qquad w_i=b\overline{z_i},\qquad
 s_{ij}=\lambda|\zeta_i-\zeta_j|^2\ne0.
\end{equation}
For a similarity $\zeta_i'=\alpha\zeta_i$,
$\lambda'=\lambda/|\alpha|^2$, one has
\[
 a'=\alpha a,\qquad z_i'=z_i,\qquad
 b'=\frac{\lambda}{|\alpha|^2}|\alpha a|^2=b,
 \qquad w_i'=w_i.
\]
Conversely,
\[
 z_i'=z_i\ (\forall i)
 \ \Longrightarrow\
 \zeta_i'=\frac{a'}a\zeta_i\ (\forall i),
\]
so the map is injective on similarity classes. Finally, any point in
the stated locus satisfies
\[
 b z_i=\overline{w_i}
 =\sum_{j\ne i}\frac{\Gamma_j}{\overline{z_i}-\overline{z_j}}
 =\sum_{j\ne i}\Gamma_j\frac{z_i-z_j}{|z_i-z_j|^2},
\]
which is \eqref{eq:physical} with positions $z_i$, center $0$ and
angular parameter $b$.
\end{proof}

\section{Isolation and the B\'ezout bound}

We first prove an isolation criterion by logarithmic residues, then bound the number of isolated solutions
by B\'ezout and use the same estimate to obtain an explicit finite union
of exceptional quadratic hypersurfaces.

Eliminate $z_2$ by $z_2=z_1+1$. The remaining coordinates are
\[
 (z_1,z_3,\ldots,z_N,w_1,\ldots,w_N)\in\C^{2N-1};
\]
we continue to abbreviate such a point by $(z,w)$, with $z_2=z_1+1$
understood in every formula. All vortex indices range from $1$ to $N$;
an empty product equals $1$. For $1\le i\le N$, define
\begin{align}
 A_i&=z_i\prod_{j\ne i}(w_i-w_j)
       -\sum_{j\ne i}\Gamma_j
          \prod_{k\notin\{i,j\}}(w_i-w_k),\notag\\
 B_i&=w_i\prod_{j\ne i}(z_i-z_j)
       -\sum_{j\ne i}\Gamma_j
          \prod_{k\notin\{i,j\}}(z_i-z_k).
 \label{eq:polynomials}
\end{align}
Clearing the nonzero denominators in \eqref{eq:complex} gives
\begin{equation}\label{eq:open-algebraic}
\begin{split}
 Y_{\G}&:=\{(z,w)\in\C^{2N-1}:A_i=B_i=0\text{ for }1\le i\le N\},\\
 X_{\G}&=\{(z,w)\in Y_{\G}:s_{ij}\ne0\text{ for every }i<j\},
 \qquad s_{ij}=(z_i-z_j)(w_i-w_j).
\end{split}
\end{equation}
An \emph{isolated solution} means an isolated point of $X_{\G}$ in the
complex Euclidean topology. Since $X_{\G}$ is an open subset of the
algebraic set $Y_{\G}$, its isolated points are precisely the
zero-dimensional irreducible components of $Y_{\G}$ that satisfy
all the inequalities $s_{ij}\ne0$. An algebraic set has finitely many
irreducible components; thus, if every solution is isolated, $X_{\G}$
is finite. All counts below concern distinct points; no Jacobian
nondegeneracy is assumed.

\subsection{Isolation by logarithmic residues}

\Needspace{10\baselineskip}
\begin{proposition}\label{prop:residue-isolation}
Fix $\G\in(\C^*)^N$. If $X_{\G}$ has a nonisolated point, then
there exists
\[
 0\ne k\in\Z^{\binom N2}
 \qquad\text{such that}\qquad
 Q_k(\G)=\sum_{i<j}k_{ij}\Gamma_i\Gamma_j=0.
\]
Consequently, if the numbers $\Gamma_i\Gamma_j$, $1\le i<j\le N$,
are linearly independent over $\Q$, every point of $X_{\G}$ is isolated.
\end{proposition}

\begin{proof}
Let $x\in X_{\G}$ be nonisolated, and let $Z$ be a positive-dimensional
irreducible component of $Y_{\G}$ containing $x$. Put $d=\dim Z\ge1$.

\smallskip
\emph{An algebraic curve through $x$.}
Let $\overline Z$ denote the projective closure of $Z$ in complex
projective space $\PP^{2N-1}(\C)$, with $\C^{2N-1}$ identified with its usual affine
coordinate chart. If $d=1$, take $C=\overline Z$. If $d\ge2$, choose
projective hyperplanes
\[
 H_1,\ldots,H_{d-1}\subset\PP^{2N-1}(\C)
\]
successively through $x$, each containing no irreducible component of
the preceding intersection. Such choices exist because every component
at that stage is positive-dimensional, and the hyperplanes through $x$
containing it form a proper closed subset of the hyperplane space.
Each intersection has pure dimension one less than the preceding one.
Choose a reduced irreducible component $C$ through $x$ of
\begin{equation}\label{eq:curve-section}
 \overline Z\cap H_1\cap\cdots\cap H_{d-1};
 \qquad \dim C=1.
\end{equation}
For $d=1$, this intersection means $\overline Z$ itself. Since
$\overline Z\cap\C^{2N-1}=Z\subset Y_{\G}$,
\[
 C\cap X_{\G}
 =\{(z,w)\in C\cap\C^{2N-1}:s_{ij}\ne0\text{ for every }i<j\}
 \ni x
\]
is a nonempty open subset of $C$, with finite complement.
Thus $C$ is algebraic by construction.

Let $\nu:\widehat C\to C$ be its normalization: $\widehat C$ is a
nonsingular projective curve, and $\nu$ has finite fibers and is an isomorphism away from
the finitely many singular points of $C$
\cite[I, Proposition~6.8 and Corollary~6.11]{Hartshorne1977}.
Over $\C$, $\widehat C$ is a compact connected Riemann surface.
We use the same symbols $z_i,w_i,s_{ij}$ for their pullbacks to
meromorphic functions on $\widehat C$.
The force equations \eqref{eq:complex} hold as identities of
meromorphic functions, and $s_{ij},z_i-z_j,w_i-w_j$ are not
identically zero for $i<j$.

We next claim that {\it there is at least one $s_{ij}$ that is non constant on $C$.}
\begin{proof}[Proof of the claim]
    
Suppose all $s_{ij}$ are constant. They are nonzero on $C\cap X_{\G}$.
Put
\[
 u_i=z_i-z_1,\quad v_i=w_i-w_1,\qquad
 u_1=v_1=0,\quad u_2=1,\quad v_2=s_{12}=:b\in\C^*.
\]
For $i\ge3$,
\begin{align*}
 u_iv_i&=s_{1i},\qquad (1-u_i)(b-v_i)=s_{2i},\\
 v_i&=s_{1i}+b-s_{2i}-bu_i,\\
 0&=bu_i^2-(s_{1i}+b-s_{2i})u_i+s_{1i}
   =b(u_i-\alpha_i)(u_i-\beta_i),
\end{align*}
where $\alpha_i,\beta_i\in\C$ are the roots of this quadratic.
The identity theorem implies that one of the two factors vanishes
identically. Hence
\[
 u_i\in\{\alpha_i,\beta_i\}\subset\C,\qquad
 v_i=s_{1i}/u_i\in\C.
\]
Thus every $u_i,v_i$ is constant, and
\[
 z_i=\sum_{j\ne i}\frac{\Gamma_j}{v_i-v_j}\in\C,
 \qquad
 w_i=\sum_{j\ne i}\frac{\Gamma_j}{u_i-u_j}\in\C.
\]
All affine coordinates of $C$ would be constant, contradicting
$\dim C=1$. Hence at least one $s_{ij}$ is nonconstant.
\end{proof}

\smallskip
\emph{The residue relation.}
Pairing $(i,j)$ with $(j,i)$ in the force equations gives
\begin{equation}\label{eq:complex-virial}
 \sum_i\Gamma_i z_iw_i
 =\sum_{i<j}\Gamma_i\Gamma_j
    \left(\frac{w_i}{w_i-w_j}+\frac{w_j}{w_j-w_i}\right)
 =\sum_{i<j}\Gamma_i\Gamma_j.
\end{equation}
Since the circulations $\Gamma_i$ are fixed, differentiation on $\widehat C$ yields
\begin{align}
 0
 &=\sum_i\Gamma_i(w_i\,dz_i+z_i\,dw_i)\notag\\
 &=\sum_i\sum_{j\ne i}\Gamma_i\Gamma_j
     \left(\frac{dz_i}{z_i-z_j}+\frac{dw_i}{w_i-w_j}\right)\notag\\
 &=\sum_{i<j}\Gamma_i\Gamma_j
     \left(\frac{dz_i-dz_j}{z_i-z_j}
                  +\frac{dw_i-dw_j}{w_i-w_j}\right)\notag\\
 &=\sum_{i<j}\Gamma_i\Gamma_j\frac{ds_{ij}}{s_{ij}}.
 \label{eq:differential}
\end{align}
Choose a nonconstant $s_{ij}$. It has a pole $p\in\widehat C$:
otherwise it would be holomorphic on the compact Riemann surface
$\widehat C$, hence constant by the maximum principle.

Choose a local holomorphic coordinate $t$ with $t(p)=0$.
For any nonzero meromorphic function $g$ near $p$, write
\[
 g=t^m h(t),\qquad m\in\Z,\qquad
 h\text{ holomorphic},\quad h(0)\ne0.
\]
The \emph{order} $\ord_p g:=m$ is positive at a zero and negative at a
pole. The \emph{residue} $\Res_p$ of a meromorphic differential is the
coefficient of $dt/t$ in its local expansion. Thus
\begin{equation}\label{eq:local-residue}
 \frac{dg}{g}=m\frac{dt}{t}+\frac{h'(t)}{h(t)}\,dt,
 \qquad \Res_p\frac{dg}{g}=m=\ord_p g.
\end{equation}
Taking the residue of \eqref{eq:differential} at this fixed point yields
\begin{equation}\label{eq:residue}
 k:=(\ord_p s_{ij})_{i<j}\in\Z^{\binom N2}\setminus\{0\},
 \qquad
 0=\sum_{i<j}\Gamma_i\Gamma_j\ord_p s_{ij}=Q_k(\G).
\end{equation}
The vector $k$ is nonzero because the chosen $s_{ij}$ has a pole at $p$.
Rational linear independence of the pair products rules out
\eqref{eq:residue}, proving the isolation assertion.
\end{proof}

\subsection{B\'ezout: the number of isolated solutions}

For an irreducible projective algebraic set $W\subset\PP^n$ of
dimension $d$, its \emph{degree} $\deg W$ is the number of points of
$W\cap L$, counted with multiplicity, where $L\subset\PP^n$ is a
general linear subspace of dimension $n-d$
\cite[Definition~18.1(iii), p.~225]{Harris1992}.
For an irreducible affine algebraic set $W\subset\C^n$, set
$\deg W:=\deg\overline W$, where $\overline W\subset\PP^n$ is its
projective closure under $x\mapsto[1:x]$. In particular, a point has
degree $1$. All sums below run over irreducible components of all
dimensions; an empty sum is zero. The degree of a polynomial is its
total degree.

\Needspace{10\baselineskip}
\begin{lemma}\label{lem:bezout}
Let $Y\subset\C^n$ be an algebraic set and $P$ a polynomial of degree
at most an integer $e\ge1$. Let $Z_1,\ldots,Z_r$ be the irreducible
components of $Y$, and $W_1,\ldots,W_s$ those of $Y\cap\{P=0\}$.
Then
\[
 \sum_{b=1}^s\deg W_b\le e\sum_{a=1}^r\deg Z_a.
\]
The same inequality holds for a projective algebraic set cut by a
homogeneous polynomial. In particular, for every $\G\in(\C^*)^N$,
if $Z_1,\ldots,Z_r$ are the irreducible components of $Y_{\G}$, then
\begin{equation}\label{eq:degree}
 \sum_{a=1}^r\deg Z_a\le N^{2N}=D_N.
\end{equation}
\end{lemma}

\begin{proof}
Consider one component $Z_a$. If $P$ vanishes identically on $Z_a$,
its degree remains unchanged and $\deg Z_a\le e\deg Z_a$.
If $P$ does not vanish identically and either $\dim Z_a=0$ or $P$ is
a nonzero constant, the intersection is empty. In the remaining case,
homogenize $P$ and apply projective B\'ezout
\cite[I, Theorem~7.7]{Hartshorne1977}. Discarding intersection
multiplicities and components at infinity gives
\[
 \sum_W\deg W\le e\deg Z_a,
\]
where $W$ runs over the irreducible components of $Z_a\cap\{P=0\}$.
Summing over $a$ proves the inequality; discarding repeated components
or components contained in others can only decrease the sum.
The same argument applies to projective intersections.
Since $\deg A_i,\deg B_i\le N$, apply the inequality successively
to the $2N$ polynomials in \eqref{eq:polynomials}, starting with
$\C^{2N-1}$, of degree $1$. This gives \eqref{eq:degree}.
\end{proof}

\Needspace{9\baselineskip}
\begin{proposition}\label{prop:isolated-count}
For every $\G\in(\C^*)^N$,
\[
 \#\{x\in X_{\G}:x\text{ is isolated in }X_{\G}\}\le D_N=N^{2N}.
\]
This counts isolated points even when $X_{\G}$ also has
positive-dimensional components. In particular, if every point of
$X_{\G}$ is isolated, then $\#X_{\G}\le D_N$.
\end{proposition}

\begin{proof}
Let $Z_1,\ldots,Z_r$ be the irreducible components of $Y_{\G}$.
Since $X_{\G}$ is open in $Y_{\G}$, each isolated point of $X_{\G}$
is a zero-dimensional component of $Y_{\G}$. Hence
\[
 \#\{x\in X_{\G}:x\text{ is isolated}\}
 \le\#\{a:\dim Z_a=0\}
 \le\sum_{a=1}^r\deg Z_a\le D_N.\qedhere
\]
\end{proof}

\subsection{Bounding the residue coefficients}

The same degree estimate bounds the integers in the residue obstruction.
This reduces the possible relations to the finite list in
\eqref{eq:exceptional} and proves finiteness outside $\E_N$.

\Needspace{7\baselineskip}
\begin{corollary}\label{prop:complex}
For every $\G\in(\C^*)^N\setminus\E_N$, every point of $X_{\G}$
is isolated and
\[
 \#X_{\G}\le D_N.
\]
\end{corollary}

\begin{proof}
Suppose $X_{\G}$ has a nonisolated point. Take the component $Z$,
the curve $C$ with normalization $\nu:\widehat C\to C$, and the
pole $p$ used to construct the residue vector $k$ in the proof of
Proposition~\ref{prop:residue-isolation}. Applying Lemma~\ref{lem:bezout}
with $e=1$ at each hyperplane intersection in \eqref{eq:curve-section}
gives
\begin{equation}\label{eq:curve-degree}
 \delta:=\deg C\le\deg\overline Z=\deg Z\le D_N.
\end{equation}

Consider a coordinate difference $f=z_i-z_j$ or $f=w_i-w_j$ with
$i<j$, as a meromorphic function on $\widehat C$.
If $f$ is constant, it is nonzero, and $\ord_p f=0$.
Otherwise choose $a\in\C$ avoiding the finite values of $f$ at the
finitely many points of $\widehat C$ mapped by $\nu$ to singular
points or points at infinity of $C$. Every zero of $f-a$ then lies
over a smooth affine point of $C$. Its order is the local intersection
multiplicity with the hyperplane $f=a$. B\'ezout therefore gives
\[
 \sum_{q\in\widehat C}\max\{\ord_q(f-a),0\}\le\delta.
\]
For any nonzero meromorphic function $g$ on $\widehat C$, the residue
theorem and \eqref{eq:local-residue} give
\[
 \sum_{q\in\widehat C}\ord_q g
 =\sum_{q\in\widehat C}\Res_q\frac{dg}{g}=0;
\]
only finitely many terms are nonzero. Since $f-a$ and $f$ have the
same poles with the same orders,
\[
\begin{aligned}
 \sum_{q\in\widehat C}\max\{\ord_q f,0\}
 =\sum_{q\in\widehat C}\max\{-\ord_q f,0\}=\sum_{q\in\widehat C}\max\{\ord_q(f-a),0\}\le\delta.
\end{aligned}
\]
Thus $|\ord_p f|\le\delta$ in both the constant and nonconstant cases.
The order of a product is the sum of the orders, so
\begin{equation}\label{eq:orders}
 |k_{ij}|=|\ord_p s_{ij}|
 =|\ord_p(z_i-z_j)+\ord_p(w_i-w_j)|
 \le2\delta\le2D_N.
\end{equation}
Together with \eqref{eq:residue}, this gives
\[
 0<\|k\|_\infty\le2D_N,\qquad Q_k(\G)=0,
\]
contrary to $\G\notin\E_N$. Thus every point is isolated;
Proposition~\ref{prop:isolated-count} gives the count.
\end{proof}

\begin{proof}[Proof of Theorem~\ref{thm:signed} and
Corollary~\ref{cor:zero-total}]
The monomials $\Gamma_i\Gamma_j$, $i<j$, are distinct, so
$Q_k\not\equiv0$ for every $k\ne0$. Since the index set in
\eqref{eq:exceptional} is finite, the polynomial
\[
 P_N(\G):=
 \prod_{\substack{k\in\Z^{\binom N2}\\0<\|k\|_\infty\le2D_N}}Q_k(\G)
\]
is nonzero and has integer coefficients. Thus
$\E_N=\{\G\in\C^N:P_N(\G)=0\}$ is a proper algebraic subset
of $\C^N$, defined over $\Q$.

On $H$, substitute $\Gamma_N=-\sum_{i<N}\Gamma_i$. Then
\begin{equation}\label{eq:restriction}
\begin{aligned}
 &Q_k\!\left(\Gamma_1,\ldots,\Gamma_{N-1},
                    -\sum_{i<N}\Gamma_i\right)=-\sum_{i<N}k_{iN}\Gamma_i^2
          +\sum_{i<j<N}(k_{ij}-k_{iN}-k_{jN})\Gamma_i\Gamma_j.
\end{aligned}
\end{equation}
If this polynomial vanished identically, the square coefficients would
first give $k_{iN}=0$ for $i<N$; the remaining mixed coefficients
would then give $k_{ij}=0$ for $i<j<N$. This contradicts $k\ne0$.
Each factor, and hence $P_N$, therefore remains nonzero after the
substitution. In particular, $H\cap\E_N$ is a proper subset of $H$.

The zero set of a nonzero real polynomial is closed, has empty interior,
and has Lebesgue measure zero. Apply this fact to $P_N$ on $\R^N$
and to the substituted polynomial in the $N-1$ linear coordinates
$\Gamma_1,\ldots,\Gamma_{N-1}$ on $H\cap\R^N$.
This proves the asserted nonemptiness, openness, density and full measure
in each $\mathcal O_\varepsilon$ and each nonempty
$H\cap\mathcal O_\varepsilon$, with relative topology and measure
on $H\cap\R^N$. The condition $P_N(\G)\ne0$ also gives the
Zariski-open intersections and relative genericity.

\Needspace{6\baselineskip}
Finally, for $\G\in(\R^*)^N\setminus\E_N$,
Lemma~\ref{lem:injection} and Corollary~\ref{prop:complex} give
\[
 \#\{\text{labeled physical similarity classes for }\G\}
 \le\#X_{\G}\le D_N=N^{2N}.
\]
The correspondence in Lemma~\ref{lem:injection} includes both signs
of $\lambda$ and is valid also when $\sum_i\Gamma_i=0$.
\end{proof}

\end{document}